\documentclass[11pt, a4paper]{article}
\usepackage{amsmath,amsthm,amssymb,amscd,epsfig,enumerate,graphicx,amsfonts,latexsym}

\usepackage{graphicx}
\usepackage{epsfig}
\usepackage{color}

\numberwithin{equation}{section}

\newtheorem{thm}{Theorem}[section]
\newtheorem{prop}[thm]{Proposition}
\newtheorem{lem}[thm]{Lemma}
\newtheorem{cor}[thm]{Corollary}
\theoremstyle{definition}
\newtheorem{dfn}[thm]{Definition}
\newtheorem{dfns}[thm]{Definitions}

\newtheorem{rmk}[thm]{Remark}

\newtheorem{Hyp}[thm]{Hypothesis}
\newtheorem{Not}[thm]{Notation}
\DeclareMathOperator{\Fix}{Fix}
\DeclareMathOperator{\ab}{ab}

\def\coloneq{\mathrel{\mathop\mathchar"303A}\mkern-1.2mu=}

\newcommand{\binomr}[2]{\genfrac{(}{)}{0pt}{}{#1}{#2}}
\newcommand{\gen}[1]{\left\langle#1\right\rangle}

\newcommand{\gp}[2]{\gen{#1\mid #2}}
\newcommand{\norm}[1]{\left\lVert#1\right\rVert}

\newcommand{\eps}{\varepsilon}
\newcommand{\N}{\mathbb{N}}
\newcommand{\Z}{\mathbb{Z}}
\newcommand{\R}{\mathbb{R}}

\def\coloneqq{\mathrel{\mathop\mathchar"303A}\mkern-1.2mu=}

\begin{document}
\author{Y. Antol\'{i}n,  L. Ciobanu and N. Viles}
\date{}
\title{On the asymptotics of visible elements and homogeneous equations in surface groups}
\maketitle
\begin{abstract}

Let $F$ be a group whose abelianization is $\Z^k$, $k\geq 2.$ An element of $F$
is called visible if its image in the abelianization is visible, that is, the greatest common
divisor of its coordinates is 1.

In this paper we compute three types of densities, annular, even and odd spherical, of visible elements in surface groups. We then use our results to show that the probability of a homogeneous
equation in a surface group to have solutions is neither 0 nor 1, as the lengths of the right-
and left-hand side of the equation go to infinity.
\medskip



\end{abstract}
\section{Introduction}\label{S:intro}

Let $F$ be a group whose abelianization is $\Z^k$, with $k\geq 2.$
An element of $F$ is called \textit{visible} with respect to a basis
of $\Z^k$ if its image in the abelianization is visible, that is,
the greatest common divisor of its coordinates is 1.
 Being visible is, in fact, independent
of the basis of $\Z^k$ (see Definition \ref{dfn:vis}), and we
therefore omit the references to the basis henceforth.

Let $\Sigma$ be a compact connected orientable surface of genus $r$, $r\geq 2$.  If
$\Sigma$ has no boundary, then a presentation for the fundamental group of $\Sigma$, which we call the {\it surface group of genus $r$}, is
$\gp{a_1,b_1,\dots,a_r,b_r}{[a_1,b_1]\cdots [a_r,b_r]}.$
If $\Sigma$ has boundary, then the fundamental group of $\Sigma$ is
simply a free group of finite rank.
For a group $G$, a positive integer $n$, and a fixed generating set $A$, one defines \textit{the sphere of radius} $n$ to be the set of elements of length $n$, with respect to $A$, in $G$. Then the \textit{spherical density} of a set $S$ of elements in $G$ measures the proportion of elements of length $n$ in $S$ in the sphere of radius $n,$ as $n$ goes to infinity (see Section 2). The \textit{annular density} of a set $S$ records the proportions of $S$ in two successive spheres.

While the spherical density of visible elements does not exist for the groups we consider, one can  instead look at the `odd spherical density' and `even spherical density' of visible elements of odd and even length, respectively.
In this paper we compute the annular, odd and even spherical densities of visible elements in a class of groups containing the surface groups of compact connected orientable surfaces, with or without boundary.
In \cite{KRSS07} the annular density of visible elements was computed for all free groups of finite rank (\cite{KRSS07}, Theorem A), and odd and even spherical density values were also given for the free group of rank two (\cite{KRSS07}, Theorem 3.7). Since the limits we obtain are different from $0$ and from $1$, this shows that visible elements form a set of \textit{intermediate} density in the groups we study. Intermediate density of sets in groups has been displayed for the first time in \cite{KRSS07}, and this tends to be a relatively rare behaviour for many combinatorial and algebraic properties encountered in group theory. Most of the properties studied in the literature (see for example \cite{KSS}) turned out to be \textit{negligible} or \textit{generic}, that is, with density equal to 0 or 1, respectively.

We would also like to mention the results of \cite{pr08}, where densities of sets of conjugacy classes in free and surface groups are investigated. More precisely, the density considered in \cite{pr08} is the asymptotic density of sets of root-free conjugacy classes of hyperbolic elements in surface groups, and for free groups, the density is similar to the annular density, but records the proportion in two successive balls instead of two succesive spheres. 

A consequence of our results is the fact that the solvability of homogeneous equations in the class of groups we study is a non-negligible and non-generic property. Let $G$ be a finitely generated group, $A$ a fixed generating set, and $X=\{X_1,\dots, X_n \}$, $n\geq 1$, a set of variables.
An \emph{equation} in variables $X_1, \dots, X_n$ with coefficients
$g_1,\dots, g_{m+1}$ in $G$ is a formal expression given by
$$g_1X_{i_1}^{\eps_1}g_2 X_{i_2}^{\eps_2}\dots X_{i_m}^{\eps_m}g_{m+1}=1,$$
where $m \geq 1$, $\eps_j \in \{1,-1\}$ for all $1 \leq j \leq m$,
and $i_j \in \{1, \dots, n\}.$ An equation is \textit{homogeneous} if the variables are on the left-hand side of
the equation and the constants are on the right-hand side of the equation:
\begin{equation}\label{HomEq}
X_{i_1}^{\eps_1} X_{i_2}^{\eps_2}\dots X_{i_m}^{\eps_m}=w,
\end{equation}
where $w \in G$. We say that the equation \eqref{HomEq} is {\it a homogeneous
equation of type $(m,|w|_A)$} or an {\it $(m,|w|_A)$-homogeneous equation}, where $|w|_A$ denotes the length of $w$ with respect to $A$.

We will be interested in the asymptotic behavior of $(m,|w|_A)$-homogeneous
equations when $G$ is a surface or a free group, and $m$ and $|w|_A$ go
to infinity. Our study of the asymptotics of homogeneous equations was motivated by two related questions: firstly, how often does a homogeneous equation in a free or surface group have solutions, and secondly, how likely is it, for two random words $u$ and $v$ in the group to have that $v$ is an endomorphic image of $u$? The second question was partly inspired by the work of Kapovich, Schupp and Shpilrain (\cite{KSS}). They show that the probability of two elements $u$ and $v$ in $F_k$ to be in the same automorphic orbit is $0$ as the lengths of $u$ and $v$ go to infinity. The following paragraph clarifies the relation between the two questions.

Suppose that $z(X_1, \dots, X_n)$ is the word in $X_1, \dots, X_n$
representing the left-hand side of \eqref{HomEq}, i.e.
$z(X_1, \dots, X_n)=X_{i_1}^{\eps_1} X_{i_2}^{\eps_2}\dots X_{i_m}^{\eps_m}$.
Let $F_n$ be the free group of rank $n$ on generators $x_1, \dots, x_n$.
Notice that the equation \eqref{HomEq} has solutions if and only if there
exists an homomorphism $\phi\colon F_n \rightarrow G$ such that
$\phi(z(x_1, \dots, x_n))=w$, where $z$ is written in the generators $x_1, \dots, x_n$.
The following ratios quantify the pairs of elements of the form $(z,w)$.

\begin{dfns}
Let $F,$ $G$ be countable groups and $l_F\colon F\to \N$ and $l_G\colon G\to \N$ be length functions, as defined in Definition 2.1. 

\begin{enumerate}
\item The {\it $(s,t)$-mapping ratio} $e_{\rho}(F,G,s,t)$ is the ratio
of the pairs of elements $(f,g)\in F\times G$ such that $l_F(f)\leq
s,$ $l_G(g)\leq t$ and with the property that $g$ is a homomorphic
image of $f,$ among all pairs $(f,g)\in F\times G$ with $l_F(f)\leq
s,$ $l_G(g)\leq t,$ that is,

$$e_{\rho}(F,G,s,t)=\dfrac{\sharp\{(f,g)\in F\times G: l_F(f)\leq s,l_G(g)\leq t,
\phi(f)=g\text{ for some }\phi\in \text{Hom}(F,G)\}}
{\sharp\{(f,g)\in F\times G: l_F(f)\leq s,l_G(g)\leq t\}}.$$

\item The {\it spherical $(s,t)$-mapping ratio} $e_{\gamma}(F,G,s,t)$ is
the ratio of the pairs of elements $(f,g)\in F\times G$ such that
$l_F(f)= s,$ $l_G(g)= t$ and with the property that $g$ is a
homomorphic image of $f$ among all pairs $(f,g)\in F\times G$ with
$l_F(f)= s,$ $l_G(g)= t,$ that is,

$$e_{\gamma}(F,G,s,t)=\dfrac{\sharp\{(f,g)\in F\times G: l_F(f)= s,l_G(g)= t,
\phi(f)=g\text{ for some }\phi\in \text{Hom}(F,G)\}}
{\sharp\{(f,g)\in F\times G: l_F(f)= s,l_G(g)= t\}}.$$
\end{enumerate}
\end{dfns}

In Section \ref{S:densitiesab} we will study the asymptotic behavior
of the $(s,t)$-mapping ratio $e_\rho(F,G,s,t)$ for $F$ and $G$
free-abelian groups with $l_G$ and $l_F$ beeing the restriction of the $||\cdot ||_p$ norm, $1\leq p\leq \infty$. We will show that the limit of $e_\rho(F,G,s,t)$, as $s$ and $t$ go to infinity is neither 0 nor 1. The
computation of the asymptotic behavior of this ratio is based on the
densities of visible elements in a free-abelian group.

In Section \ref{S:densitiesF} we study the annular, even and odd spherical densities of visible elements
in free and surface groups (Corollary \ref{cor:betas}). We obtain our main result (Theorem \ref{thm:main}) which relates
the densities of  visible points in surface and free groups with the
densities in the abelianization. 

In Section \ref{S:homeq} we study the asymptotic behavior of the spherical
$(s,t)$-mapping ratio \\
$e_{\gamma}(F,G,s,t)$ when $F$ and $G$ are free or surface
groups. We exploit the connection of $e_{\gamma}(F,G,s,t)$ with $e_\rho(F_{\ab},G_{\ab},s,t)$ to obtain upper and lower bounds of this asymptotic behavior. As a corollary, we obtain that the probability of an $(s,t)$-homogeneous
equation in a surface group to be solvable is neither 0 nor 1, as $s,t$ go to infinity (Corollary \ref{cor:homeq}).

The asymptotic behavior of equations in free abelian and free nilpotent groups is also being studied in a work by B. Gilman, A. Miasnikov and V. Romankov \cite{GMR}.

\section{Notation}

\begin{dfns}
Let $F$ be a finitely generable group, and let $A$ be a finite generating set of $F$.
If $w \in F$, then $|w|_A$ denotes the length of the shortest word in $A^{\pm 1}$ representing $w$.

For $1\leq p\leq \infty$, let $l_p\colon \Z^r\to \R$ denote the restriction to $\Z^r$ of the $|| . ||_p$-norm from $\R^r.$

A length function for a set $S$ is a function $l\colon S\to \N$ such that, for every $n\in \N$, the set $l^{-1}(\{0,1,2,\dots,n\})$ is finite. The functions $|.|_A$ and $l_p$ are examples of length functions in $F$ and $\Z^r.$
\end{dfns}

\begin{dfns}

Let $F$ be a group (or more generally, a set) and $l_F\colon F\to \N$ a length function.
\begin{enumerate}
\item Let $S \subseteq F$ and $n \geq 0$. Then
$$\rho_{l_F}(n,S)= \sharp\{x\in S : l_F(x) \leq n\},$$
and
$$\gamma_{l_F}(n,S)= \sharp\{x\in S : l_F(x) = n\}$$
denote the cardinality of the intersection of $S$ with the ball and sphere of radius $n$ in $F$, respectively.
\item Let $S \subseteq F$. The \textit{asymptotic density}
of $S$ in $F$ is
$$\bar{\rho}_{l_F}(S)=\limsup_{n \rightarrow \infty}\frac{\rho_{l_F}(n,S)}{\rho_{l_F}(n, F)}.$$
If the limit exists, then we denote it by $\rho_{l_F}(S)$ and
 we call it the \textit{strict asymptotic density}.

\item Let $S \subseteq F$. The \textit{spherical density}
of $S$ in $F$ is
$$\bar{\gamma}_{l_F}(S)=\limsup_{n \rightarrow \infty}\frac{\gamma_{l_F}(n, S)}{\gamma_{l_F}(n, F)}.$$
If the limit exists, then we denote it by $\gamma_{l_F}(S)$
and we call it the \textit{strict spherical density.}

\item  Let $S \subseteq F$. The \textit{annular density}
of $S$ in $F$ is
$$\bar{\sigma}_{l_F}(S)=\limsup_{n \rightarrow \infty}
\frac{1}{2}\left( \frac{\sharp\{x\in S : l_F(x) = n-1\}}
{\sharp \{x\in F : l_F(x) = n-1\}}+\frac{\sharp\{x\in S : l_F(x) = n\}}
{\sharp \{x\in F : l_F(x) = n\}}\right)$$

If the limit exists, then we denote it by $\sigma_{l_F}(S)$
and we call it the \textit{strict annular density.}
\end{enumerate}

\end{dfns}

When $F$ is a group, finitely generated by $A$, and $l_F=|\cdot|_A,$ the word length, we will just write $\rho_A, \gamma_A$ and $\sigma_A.$ Similarly if $F=\Z^r$ and  $l_F=l_p,$ the restriction of the $p$-norm, we will just write $\rho_p, \gamma_p$ and $\sigma_p.$

\begin{dfns}\label{dfn:vis}
For a nonzero element $z\in \Z^r$ we denote by $\gcd(z)$ the
 greatest common divisor of its coordinates. If $z=(0,\dots,0)\in
 \Z^{r}$ we set $\gcd(z)=\infty.$ Note that $\gcd$ is invariant
 under the action of $\text{Aut}(\Z^r)=SL(r,\Z).$ Hence, for all
 $z\in \Z^r$, $\gcd(z)$ does not depend on the basis of $\Z^r.$

An element of $z\in \Z^r$ is called {\it visible} if $\gcd(z)=1.$ If
 $\gcd(z)=t,$ then we call the element $t$-\textit{visible}.

We denote by $F_{\ab}$ the abelianization of the group $F,$
 that is, $F_{\ab}=F/[F,F].$
Suppose that $F_{\ab}$ is a free-abelian group of finite rank and let
 $\ab\colon F \rightarrow F_{\ab}$ be the abelianization map.
We say that an element $f\in F$ is {\it visible } (resp. $t$-visible)
  if $\ab(f)$ is visible (resp. $t$-visible) in $F_{\ab}$.
\end{dfns}

\section{Densities of visible elements in $\Z^{r}$}\label{S:densitiesab}

Let $r\geq 2$ be an integer and let $U_t$ denote the set of all $t$-visible elements in $\mathbb{Z}^{r}$.
For a complex number $k$, recall that the  Riemann zeta function is
given by $$\zeta(k)=\sum_{n=1}^{\infty} \frac{1}{n^k},\qquad
\mathfrak{Re}(k)>1.$$

A classical result in number theory provides the value for the strict asymptotic density of $t$-visible elements in $\Z^{r}$.

\begin{prop}[\cite{c56}]\label{prop:visible}
For any integer $t \geq 1$ $$\rho_{\infty}(U_t)= \frac{1}{t^r \zeta(r)}.$$\hfill\qed
\end{prop}

By \cite[Theorem A (1)]{KRSS07} or Remark \ref{rem:norms}, one can substitute $\rho_\infty$ by $\rho_p$ for the sets $U_t$:

\begin{prop}\label{prop:l1visible}\cite[Theorem A (1)]{KRSS07}
For any integer $t \geq 1$ and any $p,$ $1\leq p\leq \infty,$  $$\rho_p(U_t)= \rho_\infty(U_t).$$\hfill\qed
\end{prop}

The following lemma shows that homomorphisms between groups with free-abelian abelianization (of finite rank) send $t$-visible elements to $tm$-visible elements, where $t, m$ are positive integers. The second part of the lemma shows that a visible element in a group can be mapped to any element in the image via a homomorphism.

\begin{lem}\label{lem:gcd}
Let $F,G$ be groups whose abelianization is free-abelian of finite rank. Let $f\in F.$
\begin{enumerate}[{\normalfont(i).}]
\item\label{it:gcd1} Let $\phi:F \rightarrow G$ be a group homomorphism.
 Then $\gcd(\ab(\phi(f)))$ is a multiple of $\gcd(\ab(f))$. In particular, if $\gcd
(\ab(f))=\infty$, then $\gcd(\ab(\phi(f)))=\infty$.

\item\label{it:gcd2}  If, moreover $\gcd(\ab(f))=1,$ then for any element
$g$ in $G$ there exists an homomorphism $\phi\colon F \rightarrow G$ such
that $\phi(f)=g$.
\end{enumerate}
\end{lem}

\begin{proof}
Let $n$ be the rank of $F_{\ab}$ and let $\{e_1,\dots, e_n\}$ be a basis
of $F_{\text{ab}}.$ For $f\in F_{\text{ab}},$ we denote by $(f)_i$ the $i$th
 coordinate of $f$ with respect to the basis. That is, $f=(f)_1e_1+\dots+(f)_ne_n.$

\begin{enumerate}[{\normalfont(i).}]
\item
Let $g=\phi(f)$. Then $(\ab(g))_j=\sum_{i=1}^n (\ab(f))_i (\phi(e_i))_j$.

Thus each $(\ab(g))_j$ is a multiple of $\gcd(\ab (f))$, since each $(\ab(f))_i$
 is a multiple of $\gcd(\ab(f))$.

\item Since $\gcd(\ab(f))=1$, then $\gcd((\ab(f))_1, \dots, (\ab(f))_n)=1$
and therefore there exist integers $p_1, \dots, p_n$ such that $\sum_{i=1}^n
(\ab(f))_ip_i=1$. Consider the homomorphism $\psi_1\colon F_{\ab} \rightarrow
\gp{x}{\quad}$ which sends $e_i$ to $x^{p_i}$ for all $1\leq i \leq n$. It follows
that $\psi_1(\ab(f))=x$. Let $\psi_2\colon\gp{x}{\quad}\to G$
be any homomorphism sending $x$ to $g$. This shows that the composition of
$\ab,$ $\psi_1$ and $\psi_2$ produces a homomorphism $\phi\colon F \rightarrow G$
such that $\phi(f)=g$.
\end{enumerate}
\end{proof}

\begin{cor}\label{cor:homimage} Let $\Z^n$ and $\Z^k$ be the free-abelian  groups
 of ranks $n$ and $k$, respectively. Then the following inequalities hold with respect to  $l_p$ for $1\leq p\leq \infty$:
\begin{equation}\label{eq:liminfab}
\frac{1}{\zeta(n)} \leq
\liminf_{s \rightarrow \infty, t \rightarrow \infty} e_{\rho}(\Z^n,\Z^k,s,t),
\end{equation}

\begin{equation}\label{eq:limsupab}
\limsup_{s \rightarrow \infty, t \rightarrow \infty} e_{\rho}(\Z^n,\Z^k,s,t)
\leq  1- \frac{1}{\zeta(k)}\left(1-\frac{1}{\zeta(n)}\right).
\end{equation}

\end{cor}
\begin{proof}
We fix some $p,$ $1\leq p\leq \infty.$
Let $e_{\ab}(s,t)\coloneqq e_{\rho}(\Z^n,\Z^k,s,t)$ with respect the length $l_p$ and let $|u|=l_p(u)$.

By Lemma \ref{lem:gcd}\eqref{it:gcd2}
$$e_{\ab}(s,t)\geq \frac{\{(u,v)\in \Z^n\times \Z^k : |u|\leq s, |v|\leq t, \gcd(u)= 1\}}
{\rho_p(s,\Z^n)\rho_p(t,\Z^k)}=\frac{\{u\in \Z^n:|u|\leq s, \gcd(u)=1\}}{\rho_p(s,\Z^n)}.$$
Taking limits, we obtain \eqref{eq:liminfab} by Propositions \ref{prop:visible} and \ref{prop:l1visible}.

By Lemma \ref{lem:gcd}\eqref{it:gcd1}
\begin{align*}e_{\ab}(s,t)&
\leq  1-\frac{\{(u,v)\in \Z^n\times \Z^k :|u|\leq s, |v|\leq t, \gcd(u)\neq 1,
\gcd(v)=1\}}{\rho_p(s,\Z^n)\rho_p(t,\Z^k)}\\ &=1-\left(1-\frac{\{u\in \Z^n: |u|\leq s,
\gcd(u)=1\}}{\rho_p(s,\Z^n)}\right) \frac{\{v\in \Z^k :  |v|\leq t, \gcd(v)=1\}}
{\rho_p(t,\Z^k)}.
\end{align*}
Taking limits, we obtain \eqref{eq:limsupab} by Propositions \ref{prop:visible} and \ref{prop:l1visible}.
\end{proof}

One of the key ingredients needed to extend the previous result to the analogue for surface groups is
determining the asymptotic density of elements of even length in $\mathbb{Z}^k.$
This was done in \cite[Proposition 3.6 ]{KRSS07} for $k=2,$ and we now compute the value for a general $k$.

\begin{prop}\label{prop:keven} Let $k\geq2,$ and let $U_1^{ev}=\{z\in U_1: l_1(z)\text{ is even} \}$  denote the set of visible
elements of even length in $\mathbb{Z}^k$.
Then $$\rho_{\infty}(U_1^{ev})=\frac{2^{k-1}-1}{2^k-1}\rho_{\infty}(U_1)
=\frac{2^{k-1}-1}{(2^k-1)\zeta(k)}.$$
\end{prop}

\begin{proof}

Let $n$ be a positive integer and let $[0,n]=\{0, 1, \dots, n\}$.
For $X_1, \dots, X_k \in
\{\mathcal A, \mathcal O, \mathcal E\}$ we denote by $X_1 X_2 \dots X_k(n)$
 the number of all $z=(z_1, \dots, z_k) \in U_1$ such that $ z_i\in [0,n]$
 and  the parity of $z_i$ is $X_i$. Here $\mathcal{A}$ stands for ``any",
$\mathcal{E}$ stands for ``even" and $\mathcal{O}$ stands for ``odd".

We will use the convention  $\underbrace{ X \dots X}_\text{k times}= X^k$,
for any $X \in \{\mathcal A, \mathcal O, \mathcal E\}$ and $k \geq 1$.

Note that $X_1 X_2 \dots X_k(n)=X_{s(1)} X_{s(2)} \dots X_{s(k)}(n)$, for
any permutation $s$  of $\{1, \dots, k\}$, and that $\mathcal{E}^k(n)=0$
for any $k, n \geq 1$.

The total number of elements in $U_1$ in $[0,n]^k$ is\begin{equation}\label{eq:|U_1|}
 \mathcal A^k(n)= \sum_{i=1}^k \binomr{k}{i} \mathcal E^{k-i} \mathcal O^i(n).
\end{equation}

Let $U_1^{ev}(n)$ be the set  $U_1^{ev}\cap [0,n]^k.$ Then \begin{equation}\label{eq:|U_1ev|} |U_1^{ev}(n)|=
\sum_{i=1}^{\left[\frac{k}{2}\right]} \binomr{k}{2i} \mathcal
E^{k-2i}\mathcal O^{2i}(n).
\end{equation}

We claim that:

\begin{equation}\label{eq:EO=OO}
\mathcal E^{k-i} \mathcal O^i(n)=\mathcal O^k(n)+o(n^k) \text{ for all $1 \leq i \leq k$.}
\end{equation}

Assume first that \eqref{eq:EO=OO} holds. From \eqref{eq:|U_1ev|} and \eqref{eq:EO=OO} we get
$$|U_1^{ev}(n)|=\sum_{i=1}^{\left[\frac{k}{2}\right]} \binomr{k}{2i}\mathcal O^k(n)+o(n^k),$$
and since $\sum_{i=1}^{\left[\frac{k}{2}\right]} \binomr{k}{2i}=2^{k-1}-1$, we
get that $$|U_1^{ev}(n)|=(2^{k-1}-1)\mathcal O^k(n)+o(n^k).$$

Since $\sum_{i=1}^{k} \binomr{k}{i}=2^{k}-1,$ from \eqref{eq:|U_1|} and \eqref{eq:EO=OO} we get
$$\mathcal O^k(n)(2^k-1)=\mathcal A^k(n)+o(n^k),$$
and hence
$$|U_1^{ev}(n)|=\frac{2^{k-1}-1}{2^k-1}\mathcal A^k(n)+o(n^k).$$

Since $\rho_{\infty}(U_1)=\lim_{n\rightarrow \infty}\frac{\mathcal A^k(n)}{n^k}=\frac{1}{\zeta(k)}$, we get that
\begin{align*}
\rho_{\infty}(U_1^{ev})%
& =\limsup_{n \rightarrow \infty}\dfrac{|U_1^{ev}(n)|}{n^k}\\%
& =\lim_{n \rightarrow \infty}\dfrac{\frac{2^{k-1}-1}{2^k-1}\mathcal A^k(n)+o(n^k)}{n^k}\\%
&=\dfrac{2^{k-1}-1}{2^k-1}\rho_{\infty}(U_1)
=\dfrac{2^{k-1}-1}{(2^k-1)\zeta(k)}.
\end{align*}

This completes the proof of the proposition. We now  show \eqref{eq:EO=OO}.
Notice first that
$$\mathcal O^i\mathcal E^{k-i-1}\mathcal A(n)=
\mathcal O^i \mathcal E^{k-i}(n)+\mathcal O^{i+1}\mathcal E^{k-i-1}(n).$$
Hence it is enough to show
\begin{equation}
\label{eq:OA=2OE}
\mathcal O^{i}\mathcal E^{k-i-1}\mathcal A(n)=
2\mathcal O^{i}\mathcal E^{k-i}(n)+o(n^k) \text{ for all } 1\leq i \leq k.
\end{equation}
Let $\mu\colon \N\to \{-1,0,1\}$ denote the M\"obius function and recall that $\sum_{d|n}\mu(d)$
is equal to 1, if $n=1$ and 0 otherwise. Hence
$$\mathcal O^{i}\mathcal E^{k-i-1}\mathcal A(n)
=
\sum_{\stackrel{0\leq x_j\leq n,\, 2\nmid x_j}{j=1,\dots, i}}
\sum_{\stackrel{0\leq x_j\leq n,\, 2\mid x_j}{j=i+1,\dots, k-1}}
\sum_{0\leq x_k\leq n}
\sum_{\;d\mid \gcd(x_1,\dots, x_k)} \mu(d)$$ and

$$\mathcal O^{i}\mathcal E^{k-i}(n)=
\sum_{\stackrel{0\leq x_j\leq  n,\, 2\nmid x_j}{j=1,\dots, i}}
\sum_{\stackrel{0\leq x_j\leq n,\, 2\mid x_j}{j=i+1,\dots, k}}
\sum_{\;d\mid \gcd(x_1,\dots, x_k)} \mu(d).$$

Now we switch the order in the summation. We rearrange the terms
depending on $\;d\mid \gcd(x_1,\dots, x_k)$, writing $x_i=y_id$.
Since there is an odd coordinate, $2 \nmid d$. We obtain that

$$\mathcal O^{i}\mathcal E^{k-i-1}\mathcal A(n)
= \sum_{2\nmid d} \mu(d)
\sum_{\stackrel{0\leq y_j\leq n/d,\, 2\nmid y_j}{j=1,\dots, i}}
\sum_{\stackrel{0\leq y_{j}\leq n/d,\, 2\mid y_j}{j=i+1,\dots, k-1}}
\sum_{y_k\leq n/d} 1
$$

and
$$\mathcal O^{i}\mathcal E^{k-i}(n)
= \sum_{2\nmid d} \mu(d)
\sum_{\stackrel{0\leq y_j\leq n/d,\, 2\nmid y_j}{j=1,\dots, i}}
\sum_{\stackrel{0\leq y_{j}\leq n/d,\, 2\mid y_j}{j=i+1,\dots, k}}
 1.
$$

Hence $\mathcal O^{i}\mathcal E^{k-i-1}\mathcal A(n)-
2\mathcal O^{i}\mathcal E^{k-i}(n)$ is equal to
\begin{equation}\label{eq:error}
\sum_{2\nmid d} \mu(d)
\sum_{\stackrel{0\leq y_j\leq n/d,\, 2\nmid y_j}{j=1,\dots, i}}
\sum_{\stackrel{0\leq y_{j}\leq n/d,\, 2\mid y_j}{j=i+1,\dots, k-1}}
\left( \left[ \dfrac{n}{d} \right] -2\left[\dfrac{n}{2d}\right]\right).
\end{equation}

The term in parenthesis is either 0 or 1, and it is always 0 for $d >n.$
 Thus the asymptotic behavior of  \eqref{eq:error} is of type
\begin{align*}
O(\sum_{d\leq n}\sum_{\stackrel{0\leq y_j\leq n/d}{j=1,\dots,k-1}} 1)
\subseteq & O(\sum_{d=1}^n (n/d)^{k-1})\\
=&O(n^{k-1}\left(\frac{1}{k-2}-\frac{1}{(k-2)n^{k-2}}\right))\\
=& O(n^{k-1})\subset o(n^k)
\end{align*}
\end{proof}

\section{Densities of visible elements in surface groups}\label{S:densitiesF}

The main result of this section is an extension of \cite[Theorem A]{KRSS07} that allows us to compute densities of visible elements in free and surface groups. We need to fix some notation.

\begin{Not}\label{Not:groups}
For $k\geqslant 2$, we denote by $F_k$ the free group of rank $k$ and by $S_k$
the surface group of genus $k$.

We will work with the standard presentation for $F_k$,
$$\gp{a_1,\dots,a_k}{\quad},$$
and let  $A=\{a_1,\dots,a_k\}^{\pm 1}.$

A presentation for $S_k$ has the form
\begin{equation*}
\gp{a_1,b_1,\dots,a_k,b_k}{[a_1,b_1]\cdots [a_k,b_k]}.
\end{equation*}
In this case we let $A=\{a_1,b_1,\dots, a_k,b_k\}^{\pm 1}.$

Let $r$ denote the rank of the abelianization,
that is $r=k$ for $F_k$, and $r=2k$ for $S_k$.
\end{Not}

Our main result is based on the following local limit theorem of Sharp in \cite{sharp1}.
\begin{thm}\label{thm:hypo} (see Theorems 1, 3, 4 in \cite{sharp1}) Let $F$ be $F_k$ or $S_k$, and $A$ and $r$ be the corresponding generating set and rank of the abelianization of $F$, as in notation \ref{Not:groups}.

Let $\ab:F\to \Z^r$ be the abelianization map. Then there exists
a symmetric positive definite real matrix $D$ such that
\begin{equation}\label{eq:limhyp}
\lim_{n\to \infty}\left|(\det D)^{1/2}  n^{r/2}
\left(\frac{\gamma_A(n,\ab^{-1}(\alpha))}{\gamma_A(n,F)}+\frac{\gamma_A(n+1,\ab^{-1}(\alpha))}
{\gamma_A(n+1,F)}\right)
-\frac{2}{(2\pi)^{r/2}}e^{-\langle \alpha, D^{-1}\alpha\rangle/2n}\right|=0,
\end{equation}
uniformly in $\alpha\in \Z^{r}$.
\end{thm}
\begin{proof}
For $F=S_k$ this is exactly \cite[Theorem 4]{sharp1} with $\mathfrak{g}=r/2$.
For $F=F_k$ and $D$ the diagonal matrix with all entries equal to $\sigma^2$, one obtains exactly \cite[Theorem 1]{sharp1}.
\end{proof}

Since the proof of the main theorem of this section does not use the fact that $F$ is a free or surface group,  but only the conclusions of Theorem \ref{thm:hypo}, we will fix the following Hypothesis.
\begin{Hyp}\label{Hyp:thm}
Let $F$ be a group generated by a finite set $A$ such that $F_{\ab}\cong\Z^r$ and $D$ be a symmetric positive definite
real matrix such that the limit  \eqref{eq:limhyp} goes to zero uniformly  in $\alpha\in \Z^{r}$.
\end{Hyp}
By Theorem \ref{thm:hypo}, the free group $F_k$ and the surface group $S_k$ of Notation
\ref{Not:groups} satisfy the Hypothesis \ref{Hyp:thm}.

\begin{dfn}
Let $G_r$ be the set of all $M\in SL(r,\Z)$ such that $M=I_r$ in $SL(r,\Z/2\Z).$
 Then $G_r$ is a finite-index subgroup of $SL(r,\Z).$
\end{dfn}

\begin{dfn}
We say that a bounded open subset of $\R^r$ is \textit{nice} if its boundary is piecewise smooth. \end{dfn}

\begin{prop}\label{prop33sharp}\cite[Proposition 3.3.]{KRSS07}
Let $S\subseteq \Z^r$ be a $G_r$-invariant subset such that
$\delta=\rho_\infty(S)$ exists. Let $\Omega\subseteq \R^r$ be a nice bounded open set and for $t\in \R,$ $t>0$,
 let $$\mu_{t,S}(\Omega)\coloneq \dfrac{\sharp(S\cap t\Omega)}{t^r}.$$
Then we have
 \begin{equation}\label{eqmutS}
  \lim_{t\to\infty}  \mu_{t, S}(\Omega)=\delta\lambda(\Omega),
 \end{equation}
where $\lambda$ is the Lebesgue measure.
\end{prop}
Although \cite{KRSS07} indicates that the proof is similar to that
of \cite[Proposition 2.3]{KRSS07}, we include here a proof for
Proposition \ref{prop33sharp} for the sake of completeness.
\begin{proof}
Each $\mu_{t,S}$ can be regarded as a measure on $\R^r$. We prove
the result by showing that $\mu_{t,S}$ weakly converge to
$\delta\lambda$ as $t\to \infty$.

 By Helly's
theorem (see, for instance, \cite[Thm 25.9]{billingsley}),
there exists a sequence  $\{t_i\}$ with
$\lim_{i\to\infty} t_i=\infty$ such that the sequence $\mu_{t_1,
S}$, $\mu_{t_2, S}$, ... is weakly convergent to some limiting
measure. We now identify this measure by showing that for every
convergent subsequence of $\mu_{t_i, S}$ the limiting measure is
equal to $\delta\lambda$.

Indeed, we assume that $\eta=\{t_i\}$ is a sequence with
$\lim_{i\to\infty}t_i=\infty$ such that the sequence $\mu_{t_i, S}$ converges to the limiting
measure $\mu_{\eta}=\lim_{i\to\infty}\mu_{t_i,S}$. Every $\mu_{t_i, S}$ is invariant with
respect to the $G_r$-action on $\R^r$. Therefore the limiting measure $\mu_\eta$ is also
$G_r$-invariant. Moreover, the measures $\mu_{t,S}$ are dominated by the measures $\lambda_t$
defined as $\lambda_t(\Omega)=\frac{\sharp(\Z^r\cap t\Omega)}{t^r}$.

It is well known that if $\Omega \subseteq \R^r$ is a nice bounded open set, then the measures
$\lambda_t$ converge to the Lebesgue measure $\lambda$. It follows that $\mu_{\eta}$ is
absolutely continuous with respect to $\lambda$. It is also known that the natural action of
$G_r$ on $\R^r$ is ergodic with respect to $\lambda$ (see \cite{zimmer} for the proof of
 ergodicity). Therefore $\mu_\eta$ is a constant multiple $c\lambda$ of $\lambda$. The
constant $c$ can be computed for a set such as the open unit ball $B$ in the $\|\cdot\|_{\infty}$
norm on $\R^r$ defining the length function $l_\infty$ on $\Z^r$. By assumption we know that
$$\rho_{\infty}(S)=\lim_{t\to\infty}\frac{\sharp \{z\in \Z^r:z\in S\cap tB\}}{\sharp \{z\in \Z^r:z\in tB\}}=\delta.$$

We also have
$$ \lim_{t\to\infty}\frac{\sharp \{z\in \Z^r:z\in tB\}}{t^r}=\lambda(B)$$
and hence
$$ \lim_{t\to\infty}\frac{\sharp \{z\in \Z^r:z\in tB\}}{t^r}\frac{\sharp \{z\in \Z^r:z\in S\cap tB\}}{\sharp \{z\in \Z^r:z\in tB\}}=\lim_{t\to\infty}\frac{\sharp \{z\in \Z^r:z\in S\cap tB\}}{t^r}=\delta\lambda(B).$$

Therefore $c=\delta$ and $\mu_{\eta}=\delta\lambda$. The above argument shows in fact
that every convergent subsequence of $\mu_{t,S}$
converges to $\delta\lambda$ and $\lim_{t\to \infty}\mu_{t,S}=\delta \lambda$.
\end{proof}

\begin{rmk}\label{rem:norms}(see \cite[Theorem A]{KRSS07}) Let $1\leq p\leq \infty.$ The sets $U_q$ of $q$-visible elements in $\Z^r$ are $G_r$-invariant and $$\rho_p(U_q)=\rho_\infty(U_q).$$ 
\begin{proof}
Let $\Omega$ be an $l_p$ ball of radius 1.
It is well known that
$$\lambda(\Omega)=\lim_{t\to \infty}\dfrac{\sharp (\Z^r\cap t\Omega)}{t^r} $$

Then
\begin{align*}\rho_p(U_q) &= \lim_{n\to \infty}\dfrac{\sharp \{ x\in U_q: l_p(x)\leq n\}}{\sharp \{x\in \Z^r: l_p(x)\leq n\}}\\
& =\lim_{n\to \infty}\dfrac{\sharp \{ x\in U_q: l_p(x)\leq n\}}{\sharp \{x\in \Z^r: l_p(x)\leq n\}}\cdot \lim_{t\to \infty}\dfrac{\sharp (\Z^r\cap t\Omega)}{\lambda(\Omega)t^r}\\
& = \lim_{t\to \infty}\dfrac{\sharp \{ x\in U_q: l_p(x)\leq t\}}{\lambda(\Omega)t^r}\\
& = \lim_{t\to \infty}\dfrac{\sharp ( U_q \cap t\Omega)}{\lambda(\Omega)t^r}\\
& =  \dfrac{\delta \lambda(\Omega)}{\lambda(\Omega)}\\
& =  \rho_\infty(U_q).
 \end{align*}
 \end{proof}
\end{rmk}

\begin{dfn}\label{def:pn} Let $F$ be a group generated by the finite set $A$ such
 that $F_{\ab}\cong \Z^r.$

For an integer $n\geqslant 1$ and a point $x\in \R^{r}$, let $p_n$ be
given by
\begin{equation}\label{defpn}
p_n(x)=\frac12\left(\frac{\gamma_{A}(n-1,\{g\in F: \ab(g)=x\sqrt{n}\})}
{\gamma_A(n-1,F)}+\frac{\gamma_{A}(n,\{g\in F: \ab(g)=x\sqrt{n}\})}{\gamma_A(n,F)}\right).
\end{equation}
This is a distribution supported on finitely many points of
$\frac{1}{\sqrt{n}}\; \Z^{r}$.
\end{dfn}

We need the following results from \cite{rivin, sharp1} about the sequence of distributions $p_n$.

In our context, we need to restate our Hypothesis \ref{Hyp:thm}

\begin{prop}\label{prop27v2}(\cite{rivin, sharp1,KRSS07})
Let $F, A,r$ satisfy Hypothesis \ref{Hyp:thm}. Then there exists a
normal distribution $\mathfrak{N}$ with density $\mathfrak{n}$ such
that:
\begin{enumerate}[(a)]
 \item The sequence of distributions $p_n$ converges weakly to  $\mathfrak{n}$ and we have

      \begin{equation}\label{prop27eq1}
    \sup_{x\in \Z^{r}/\sqrt{n}}|n^{r/2}p_n(x)-\mathfrak{n}(x)|\longrightarrow 0, \; as \; n\to \infty.
      \end{equation}
 \item For $c>0,$ let $\overline{\Omega_c}\coloneqq\{x\in \R^{r} : \|x\|\geqslant c\}.$ Then
      \begin{equation}\label{prop27eq2}
    \lim_{c\to \infty}(\lim_{n\to\infty }\sum_{x\in \overline{\Omega_c}} p_n(x)) = 0.
      \end{equation}
\end{enumerate}
\end{prop}
\begin{proof}
Let $D$ be the matrix of Hypothesis \ref{Hyp:thm}, and let
$\mathfrak{n}(x)=\dfrac{e^{-\langle x,
{D^{-1}x}\rangle/2}}{(2\pi)^r(\det D)^{1/2}},$ the density of a
normal distribution $\mathfrak{N}.$ Firstly, we prove the limit in
(\ref{prop27eq1}).

After performing some easy computations,
\begin{align*}
|n^{r/2}&p_n(x)-\mathfrak{n}(x)|\\
&= \frac{1}{2(\det D)^{1/2}}\Big| (\det D)^{1/2}n^{r/2} \\
&
\qquad\qquad\qquad\qquad\times\left(\frac{\gamma_{A}(n-1,\ab^{-1}(x\sqrt{n}))}{\gamma_A(n-1,F)}
+\frac{\gamma_{A}(n,\ab^{-1}(x\sqrt{n}))}{\gamma_A(n,F)}\right)-
\frac{2}{(2\pi)^r}e^{-\langle \frac{\alpha}{\sqrt{n}},
\frac{D^{-1}\alpha}
{\sqrt{n}}\rangle/2}\Big|\\
&= \frac{1}{2(\det D)^{1/2}}\Big| (\det D)^{1/2}n^{r/2} \\
& \qquad\qquad\qquad\qquad\times
\left(\frac{\gamma_A(n-1,\ab^{-1}(\alpha))}
{\gamma_A(n-1,F)}+\frac{\gamma_A(n,\ab^{-1}(\alpha))}{\gamma_A(n,F)}\right)
-\frac{2}{(2\pi)^r}e^{-\langle \frac{\alpha}{\sqrt{n}},
\frac{D^{-1}\alpha}
{\sqrt{n}}\rangle/2}\Big|\\
&= \frac{1}{2(\det D)^{1/2}}\Big| (\det D)^{1/2} n^{r/2}\\
& \qquad\qquad\qquad\qquad\times
\left(\frac{\gamma_A(n-1,\ab^{-1}(\alpha))}
{\gamma_A(n-1,F)}+\frac{\gamma_A(n,\ab^{-1}(\alpha))}{\gamma_A(n,F)}\right)
-\frac{2}{(2\pi)^{r}}e^{-\langle \alpha,
D^{-1}\alpha\rangle/2n}\Big|,
\end{align*}
using the limit \eqref{eq:limhyp} of the Hypothesis
\ref{Hyp:thm}, and the fact that this limit is uniform in
$\alpha=x\sqrt{n}$, we obtain the desired result.

In order to show that the sequence of probability distributions
$\{p_n\}$ converges weakly to $\mathfrak{n}$, we use
\cite[Thm 25.8]{billingsley}, that is, it is necessary and sufficient
that for every bounded continuous function $f(x)$ on $\R^r$
\begin{equation}\label{eqlevyconthm}
 \lim_{n\to \infty}\int_{\mathbb{R}^r} f(x)p_n(x)d\lambda(x)=\int_{\mathbb{R}^r}f(x) \mathfrak{n}(x)d\lambda(x).
\end{equation}
We can write:
\begin{align*}
 \Big|\int_{\mathbb{R}^r} f(x)p_n(x)d\lambda(x)-\int_{\mathbb{R}^r}f(x) \mathfrak{n}(x)d\lambda(x)\Big|&=  \Big|\int_{\mathbb{R}^r}
 f(x)(p_n(x)-\mathfrak{n}(x))d\lambda(x)\Big|
 \\&\leqslant \int_{\mathbb{R}^r}
 |f(x)||n^{r/2}p_n(x)-\mathfrak{n}(x)|d\lambda(x).
\end{align*}

Given that $f$ is a bounded continuous function and
by the limit (\ref{prop27eq1}) proved above, the
hypothesis of the Dominated Convergence Theorem is satisfied. Applying this last
result, we obtain that
\begin{equation*}
 \Big|\int_{\mathbb{R}^r} f(x)p_n(x)d\lambda(x)-\int_{\mathbb{R}^r}f(x)
 \mathfrak{n}(x)d\lambda(x)\Big|\overset{n\to\infty}{\longrightarrow} 0,
 \end{equation*}
and the weak convergence of the sequence $\{p_n\}$ is proved.

We now prove $(b)$.
For $c>0$, let $\Omega_c=\{x\in\R^r : \|x\|<c\}$, and denote by $\overline{\Omega_c}$ the
complement of  $\Omega_c$. Then, by the weak convergence of the $p_n$ to $\mathfrak{n}$, we have that
\begin{align*}
\lim_{c\to\infty}\left(\lim_{n\to \infty}\sum_{x\in \overline{{\Omega_c}}}
p_n(x)\right)&=\lim_{c\to \infty}\left(1-\lim_{n\to\infty}\sum_{x\in \Omega_c} p_n(x)\right)
\\&=1-\lim_{c\to\infty}\int_{x\in \Omega_c}\mathfrak{n}(x)d\lambda(x)=0.
\end{align*}
\end{proof}

\begin{thm}\label{T:thm34}
 Let $\Omega\subseteq \R^r$ be a nice bounded open set. Let $S\subseteq \Z^r$
 be a $G_r$-invariant subset such that $\delta=\rho_\infty(S)$ exists.
Then there exists a normal distribution $\mathfrak{N}$ such that

 \begin{equation*}
  \lim_{n\to\infty}  \sum_{x\in S\cap \sqrt{n}\Omega}p_n(x/\sqrt{n})=\delta\mathfrak{N}(\Omega).
 \end{equation*}
\end{thm}

\begin{proof}
Note that the proof is the same as that of Theorem 3.4 in \cite{KRSS07}.
 The only difference lies in the use of Proposition \ref{prop33sharp}.

There exists a normal $\mathfrak{N}$
distribution  with density $\mathfrak{n}$ satisfying the conclusions of Proposition \ref{prop27v2}.

We have
\begin{align*}
 \sum_{x\in S\cap\sqrt{n}\Omega} p_n(x/\sqrt{n})&=
\sum_{y\in \frac{1}{\sqrt{n}}S\cap\Omega} p_n(y)
 \\&=n^{-r/2}\sum_{y\in \frac{1}{\sqrt{n}}S\cap\Omega}\mathfrak{n}(y)
\\&+n^{-r/2}\sum_{y\in \frac{1}{\sqrt{n}}S\cap\Omega}(n^{r/2}p_n(y)-\mathfrak{n}(y)).
\end{align*}
 The local limit theorem of Proposition \ref{prop27v2}(a) tells us that, 
 as $n\to \infty$, each summand
$n^{-r/2}p_n(y)-\mathfrak{n}(y)$ of the sum in
the last line above converges to zero, and hence so does their Cesaro mean.

Using the following convergence of the measures defined in
Proposition \ref{prop33sharp}, $$\lim_{n\to \infty}
\mu_{\sqrt{n},S}(\Omega)=\delta\lambda(\Omega),$$ (recall that
$\mu_{\sqrt{n},S}(\Omega)\coloneq \dfrac{\sharp(S\cap
\sqrt{n}\Omega)}{\sqrt{n}^r}$),
 we have that
$$
\lim_{n\to\infty}\sum_{x\in
{S\cap\sqrt{n}\Omega}}\frac{1}{(\sqrt{n})^{\;r}}\mathfrak{n}(x/\sqrt{n})=\int_{\Omega}\mathfrak{n}(y)
\delta d\lambda(y)=\delta\mathfrak{N}(\Omega). $$
\end{proof}

We obtain the main result of this section by basically following \cite[Theorem A]{KRSS07}. Our theorem provides the formula for the `spherical densities' of visible elements in groups that satisfy Hypothesis \ref{Hyp:thm}, which include free groups of all finite ranks and surface groups.

\begin{thm}\label{thm:main}(see also \cite[Theorem A]{KRSS07})

Let $F, A,r$ satisfy Hypothesis \ref{Hyp:thm}, 
 $S\subseteq  \Z^r$ be a $G_r$-invariant subset and $\tilde{S}=\ab^{-1}(S).$

\begin{enumerate}[{\normalfont(i).}]
\item The strict annular density $\sigma_A(\tilde{S})$ exists and, moreover,
 $\sigma_A(\tilde{S})=\rho_{\infty}(S)$.
\item Let $U_1$ denote the set of visible elements in $\Z^r$ and $V_1=\ab^{-1}(U_1)$
denote the visible elements in $F.$ Let $U_1^{ev}=\{z\in U_1: l_1(z)\text{ is even} \}$  denote the visible
elements of even length. If $\ab^{-1}(U_1^{ev})=\{v\in V_1: |v|_A \text{ is even}\}$, then

\begin{tabular}{ r c c c l }
\(\lim_{m\to\infty}\dfrac{\gamma_A(2m,V_1)}{\gamma_A(2m,F)}\)  & \(=\) &
\(2 \rho_{\infty}(U_1^{ev}) \)& \(=\) &
\(\dfrac{2^{r}-2}{(2^r-1)\zeta(r)},\)\\
 \(\lim_{m\to\infty}\dfrac{\gamma_A(2m-1,V_1)}{\gamma_A(2m-1,F)}\) &\(=\)&
\(2 \rho_{\infty}(U_1)-{2}\rho_{\infty}(U_1^{ev}) \) &\(=\)&
 \(\dfrac{2^r}{(2^r-1)\zeta(r)}.\)
\end{tabular}

\end{enumerate}
\end{thm}
\begin{proof}
For $c>0$ let $\Omega_c\coloneqq \{x\in \R^r: \norm{x}<c\}$ and let
$\overline{\Omega_c}$ be the complement of $\Omega_c.$ Then
\begin{equation}\label{eq:limnormal}
\lim_{c\to \infty}\mathfrak{N}(\Omega_c)=1
\end{equation} Let $\epsilon >0$ be arbitrary. By \eqref{eq:limnormal} and
Proposition \ref{prop27v2} (b) 
we can choose $c>0$ such that
$$|\mathfrak{N}(\Omega_c)-1|\leq \epsilon/3$$
and  
$$\lim_{n\to \infty}\sum_{x\in \overline{\Omega_c}} p_n(x)\leq \epsilon/6. $$

Let $S$ be a $G_r$-invariant subset of $\Z^r.$
By Theorem \ref{T:thm34} and the above formula there is some $n_0\geq 1$
such that  for all $n\geq n_0$ we have
\begin{equation}\label{eq:eps1}\left|\sum_{x\in S\cap\sqrt{n}\Omega_c}p_n(x/\sqrt{n})-
\rho_{\infty}(S)\mathfrak{N}(\Omega_c) \right|\leq \epsilon/3, \end{equation}
and \begin{equation}\label{eq:eps2}\sum_{x\in \overline{\Omega_c}} p_n(x)\leq \epsilon/3.\end{equation}

Let
$$Q(n)\coloneq \frac{\gamma_A(n-1, \ab^{-1}(S))}{2\gamma_A(n-1,F)}+
\frac{\gamma_A(n, \ab^{-1}(S))}{2\gamma_A(n,F)}.$$

For  $n\geq n_0$ we let 
\begin{align*}
Q(n)=&
\\
& \left(\dfrac{\sharp\{g\in F:\ab(g)\in S,|g|_A=n-1
\text{ and }\norm{\ab({g})}<c\sqrt{n}\}}{2\gamma_A(n-1,F)}  \right.\\
+& \left.\dfrac{\sharp\{g\in F:\ab(g)\in S,|g|_A=n
\text{ and }\norm{\ab({g})}<c\sqrt{n}\}}{2\gamma_A(n-1,F)}\right)\\
 +& \left(\dfrac{\sharp\{g\in F:\ab(g)\in S,|g|_A=n-1
\text{ and }\norm{\ab({g})}\geq c\sqrt{n}\}}{2\gamma_A(n-1,F)} \right.\\
+& \left.\dfrac{\sharp\{g\in F:\ab(g)\in S,|g|_A=n
\text{ and }\norm{\ab({g})}\geq c\sqrt{n}\}}{2\gamma_A(n-1,F)}\right) \\
=&\sum_{x\in S\cap \sqrt{n}\Omega_c}p_n(x/\sqrt{n})+ \sum_{x\in
S\cap (\R^r\setminus\sqrt{n}\Omega_c)} p_n(x/\sqrt{n}).
\end{align*}

In the last line of the above equation, by \eqref{eq:eps1}, the
first sum differs from $\rho_{\infty}(S)\mathfrak{N}{(\Omega_c)}$ by
at most $\epsilon/3$ since $n\geq n_0$, and  by \eqref{eq:eps2}, the
second sum is $\leq \epsilon/3$ given the choice of $c$ and $n_0$.

Therefore, again by the choice of $c$, we have $|Q(n)-\rho_{\infty}(S)|\leq \epsilon.$
Since $\epsilon$ is arbitrary, this proves (i).

We now prove (ii). First notice that since $U_1$ is $SL(r,\Z)$-invariant, it is also
$G_r$-invariant. We check that $U_1^{ev}$ is $G_r$-invariant as well. Let $u\in \Z$.
 Then $u\in U_1^{ev}$ if and only if $\sum_{1\leq i\leq r} (u)_i \mod 2=0$ and
$\gcd(u)=1.$ Let $M\in G_r.$ As $M\in SL(r,\Z),$ $\gcd(Mu)=\gcd(u)=1.$ Also, as
 $M=I_r$ in $SL(r,\Z/2\Z),$ $$\sum_{1\leq i\leq r} (Mu)_i \mod 2=\sum_{1\leq i\leq r}
 (u)_i \mod 2=0.$$   Hence, $U_1^{ev}$ is $G_r$-invariant.

We now take $S= U_1^{ev}$, for $n\geq 2$ even. Then
$$Q(n)= \dfrac{\gamma_A(n-1, \ab^{-1}(U_1^{ev}))}{2\gamma_A(n-1,F)}
+\frac{\gamma_A(n, \ab^{-1}( U_1^{ev}) )}{2\gamma_A(n,F)}$$
$$=\dfrac{\gamma_A(n, \ab^{-1}(U_1^{ev}))}{2\gamma_A(n,F)}.$$
The latter equality follows from the fact that  $\ab^{-1}(U_1^{ev})
=\{v\in V_1: |v|_A \text{ is even}\}.$

By (i), $$\lim_{m\to\infty}\frac{\gamma_A(2m,V_1)}
{\gamma_A(2m,F)}=2\lim_{m\to \infty}Q(2m)=2\rho_{\infty}(U_1^{ev}).$$
Thus $\lim_{m\to\infty}\dfrac{\gamma_A(2m-1,V_1)}{\gamma_A(2m-1,F)}=2
\rho_{\infty}(U_1)-2\rho_{\infty}(U_1^{ev}).$ 
By Proposition \ref{prop:visible} and Proposition \ref{prop:keven},
 we obtain the desired results.
\end{proof}

We now focus on surface and free groups.

\begin{cor}\label{cor:betas}Let $k\geq 2$ and let $F$ be a free group of rank
$k$ or a surface group of genus $k.$ Let $A$ and $r$ be as in Notation \ref{Not:groups}. Then
\begin{enumerate}[{\normalfont(i).}]
\item $\lim_{m\to\infty}\dfrac{\gamma_A(2m,V_1)}{\gamma_A(2m,F)}=\dfrac{2^{r}-2}{(2^r-1)\zeta(r)} .$
\item $\lim_{m\to\infty}\dfrac{\gamma_A(2m-1,V_1)}{\gamma_A(2m-1,F)}=\dfrac{2^r}{(2^r-1)\zeta(r)}.$
\end{enumerate}
\end{cor}
\begin{proof}
By Theorem \ref{Hyp:thm}, $F,A$ and $r$ satisfy the Hypothesis of Theorem \ref{thm:main}.
It only remains to show that  $\ab^{-1}(U_1^{ev})=\{v\in V_1: |v|_A \text{ is even}\}.$
Let $f$ be an element of $F$ such that $\ab(f)=0\in \Z^r.$ Then any word representing $w$
has the same number of $a$ and $a^{-1}$ and thus it has even length.

Since $\ab$ maps elements of $A$ to unit vectors, for $u\in U_1^{ev}$ there exists
$v\in \ab^{-1}(U_1^{ev})$ of even length. If $\ab(v)=\ab(v')$, then $\ab(v'v^{-1})=0$. Hence
$v'v^{-1}$ has even length, and so does $v'$.
Thus Theorem \ref{thm:main}(ii) applies.
\end{proof}

\section{Asymptotic behavior of homogeneous equations in surface groups}\label{S:homeq}

We now study the asymptotic behavior of $e_{\gamma}(G_n,G_k,s,t)$ when $G_n$ and $G_k$
 are surface or free groups, or more generally, satisfy the hypothesis of Theorem \ref{thm:main} (ii).

\begin{thm}\label{thm:homimage} Let $G_k$ and $G_n$ be free or surface groups and let $A,$ $B$
 be their respective generating sets, as in Notation \ref{Not:groups}. Let $r(k)$ and $r(n)$ denote
the ranks of the abelianization of $G_k$ and $G_n$, respectively.  Let $\epsilon,\delta\in \{0,1\}$. Then the
following inequalities hold:
$$\frac{2^{r(n)}-2(1-\eps)}{(2^{r(n)}-1)\zeta(r(n))} \leq
\liminf_{s \rightarrow \infty, t \rightarrow \infty} e_\gamma(G_n,G_k,2s+\eps,2t+\delta), $$

$$\limsup_{s \rightarrow \infty, t \rightarrow \infty} e_\gamma(G_n,G_k,2s+\eps,2t+\delta)
\leq  1- \frac{2^{r(k)}-2(1-\delta)}{(2^{r(k)}-1)\zeta(r(k))}\left(1-\frac{2^{r(n)}-2(1-\eps)}
{(2^{r(n)}-1)\zeta(r(n))}\right).$$

\end{thm}

\begin{proof}
Let $V_t$ and $W_t$ denote  the sets of  $t$-visible elements in $G_n$ and $G_k$, respectively.
Let $$E(s,t)=\{(u,v)\in G_n \times G_k : |u|_A=s, |v|_B=t, \phi(u)=v \ \textrm{for some} \
\phi \in \mathop{\text{Hom}}(G_n,G_k)\} .$$
Then $e_{\gamma}(G_n,G_k,s,t)=\frac{|E(s,t)|}{\gamma_B(s,G_n)\gamma_A(t,G_k)}.$

By Lemma \ref{lem:gcd} we have the following inequalities:

$$\gamma_B(s,W_1)\gamma_A(t,G_k) \leq  |E(s,t)|\leq \gamma_B(s,G_n)\gamma_A(t,G_k) -
\sum_{r\neq 1} \gamma_B(s,W_r)\gamma_A(t,V_1).$$

The left inequality holds because every element $v$ in $G_k$ is the homomorphic image
of a visible element in $G_n$. The right inequality holds because no visible element in
$G_k$ is the homomorphic image of an $r$-visible element in $G_n$, if $r\neq 1$.

By dividing both sides by $\gamma_B(s,G_n)\gamma_A(t,G_k)$, we get
$$\frac{\gamma_B(s,W_1)}{\gamma_B(s,G_n)}\leq e_\gamma(G_n,G_k,s,t)\leq 1-\frac{\sum_{r\neq 1}
\gamma_B(s,W_r)\gamma_A(t,V_1)}{\gamma_B(s,G_n)\gamma_A(t,G_k)}=f(s,t),$$
where
$$f(s,t)= 1 - \frac{\gamma_A(t,V_1)}{\gamma_A(t,G_k)}\frac{\gamma_B(s,G_n)-\gamma_B(s,W_1)}
{\gamma_B(s,G_n)}.$$ Let us use $\beta_{m,k}$ to denote the limits, which depend on the parity
of $m$ and the rank of the abelianization of $G_n$ and $G_k$, found in Corollary \ref{cor:betas}.
That is, $\beta_{m,k}=\frac{2^{r(k)}-2}{(2^{r(k)}-1)\zeta(r(k))}$ if $m$ is even, and $\beta_{m, k}=
\frac{2^{r(k)}}{(2^{r(k)}-1)\zeta(r(k))}$ if $m$ is odd. In order to simplify the exposition we will abuse the fact that $\beta_{m, k}$ depends on the parity of $m$ and for the next paragraph ignore
the parities of $s$ and $t$.

Then $$\lim_{s \rightarrow \infty, t \rightarrow \infty}f(s,t)=1- \beta_{t,k}(1-\beta_{s,n}),$$ and
we get the following inequalities

\begin{equation}\label{beta}\beta_{s,n} \leq  \liminf_{s \rightarrow \infty, t \rightarrow \infty}
 e_\gamma(G_n,G_k,s,t) \leq \limsup_{s \rightarrow \infty, t \rightarrow \infty} e_\gamma(G_n,G_k,s,t)\leq 1-
\beta_{t,k}(1-\beta_{s,n}).
\end{equation}

Now taking into account the parities of $s$ and $t$ we get the inequalities in the statement of the theorem.
\end{proof}
Thus the probability of an $(s,t)$-homogeneous equation to be solvable is neither 0 nor 1 as $s,t$ go to infinity. One sees this by choosing $G_n$ to be the free group on $n$ generators and $G_k$ a surface group of genus $g\geq 2$ or a free group of rank $\geq 2$ in Theorem \ref{thm:homimage}.

\begin{cor}\label{cor:homeq}
Let $G$ be a surface group of genus $g\geq 2$ or a free group of rank $\geq 2.$

Let $$A(s,t)=\dfrac{\sharp\{\text{solvable $(s,t)$-homogeneous equations in $G$ in $n$ variables}\}}
{\sharp\{\text{$(s,t)$-homogeneous equations in $G$ in $n$ variables}\}}.$$ Then
$$0<\liminf_{s \rightarrow \infty, t \rightarrow \infty} A(s,t)\leq
 \limsup_{s \rightarrow \infty, t \rightarrow \infty}A(s,t)<1.$$
\end{cor}

Similarly, by choosing both $G_n$ and $G_k$ in Theorem \ref{thm:homimage} to be surface groups one obtains the following. 
\begin{cor}\label{cor:curves}
Let $\Sigma$ be an orientable closed surface of genus $k\geq 2.$  We fix a presentation  for
$\pi_1(\Sigma),$ $\gp{a_1,b_1,\dots,a_k,b_k}{[a_1,b_1]\cdots [a_k,b_k]}.$ For a closed curve $\gamma$
in $\Sigma$ we denote by  $[\gamma]$ the image of $\gamma$ in $\pi_1(S)$ and by $|[\gamma]|$ the
length of $[\gamma]$ with respect to $\{a_1,b_1,\dots,a_k,b_k\}.$

We say that $\gamma_2$ is the image of $\gamma_1,$ if it is the image of $\gamma_1$  under a
continuous map $S\to S.$

Let
$$B(s,t)=\dfrac{\sharp\{\text{$([\gamma_1],[\gamma_2])\in \pi_1(S)^2,$
$(|[\gamma_1]|,|[\gamma_2]|)=(s,t)$
with $\gamma_2$ the image of
$\gamma_1$}\}}{\sharp\{\text{$([\gamma_1],[\gamma_2])\in \pi_1(S)^2,$
$(|[\gamma_1]|,|[\gamma_2]|)=(s,t)$}\}}.$$ Then
$$0<\liminf_{s \rightarrow \infty, t \rightarrow \infty} B(s,t)\leq
 \limsup_{s \rightarrow \infty, t \rightarrow \infty} B(s,t)<1.$$
\end{cor}
Thus for a fixed orientable surface $\Sigma$, the probability of a closed curve in $\Sigma$ to
be the image of another closed curve in $\Sigma$ by a continuous map is neither 0 nor 1, as the
curves get more and more ``complicated.''

\section*{Acknowledgments}
We are grateful to Fernando Chamizo for helpful conversations.
\medskip

\footnotesize

The first-named author was jointly funded by the MEC (Spain) and the EFRD (EU)
through Projects MTM2006-13544 and MTM2008-01550.

The second-named author was partially supported by the SNF (Switzerland) through project
number 200020-113199 and by the Marie Curie Reintegration Grant 230889.

The third-named author was supported by the grant MEC-FEDER Ref. MTM2009-08869 from the Dirección General de Investigación, MEC (Spain).

\textsc{Y. Antolin, Departament de  Matem\`atiques,
Universitat Aut\`onoma de Barcelona,
E-08193 Bella\-terra (Barcelona), Spain}

\emph{E-mail address}{:\;\;}\texttt{yagoap@mat.uab.cat}

\medskip

\textsc{L. Ciobanu,
Mathematics Department,
University of Fribourg,
Chemin du Muse\'e 23,
CH-1700 Fribourg, Switzerland
}

\emph{E-mail address}{:\;\;}\texttt{laura.ciobanu@unifr.ch}
\medskip

\textsc{N. Viles, Departament de  Matem\`atiques,
Universitat Aut\`onoma de Barcelona,
E-08193 Bella\-terra (Barcelona), Spain}

\emph{E-mail address}{:\;\;}\texttt{nviles@mat.uab.cat}

\end{document}